\documentclass[11pt]{amsart}
\usepackage{amsmath}
\usepackage{amsfonts}
\usepackage{amssymb}
\newtheorem{theorem}{Theorem}

\newtheorem{corollary}[theorem]{Corollary}

\newtheorem{remark}[theorem]{Remark}
\newtheorem{lemma}[theorem]{Lemma}
\newtheorem{prop}[theorem]{Proposition}
\newtheorem{defn}[theorem]{Definition}
\newcommand{\bb}[1]{\mathbb{#1}}
\newcommand{\cl}[1]{\mathcal{#1}}

\begin{document}

\title[On the ranges of bimodule projections]{On the ranges of
bimodule projections}
\author[A. Katavolos]{Aristides Katavolos}
\address{Department of Mathematics, University of Athens, Athens,
GREECE}
\email{akatavol@math.uoa.gr}

\author[V.~I.~Paulsen]{Vern I.~Paulsen}
\address{Department of Mathematics, University of Houston,
Houston,
Texas
77204-3476,
U.S.A.}
\email{vern@math.uh.edu}

\thanks{Research supported in part by a grant from the
University of Athens (A. Katavolos) and the NSF (V.I. Paulsen)\\
Preliminary versions of this work were announced by the authors
in
2001 and 2002, in particular during the Workshop on
Nonself-adjoint Operator Algebras held at the Fields Institute
(see \:
\texttt{http://www.fields.utoronto.ca/audio/02-03/nsoa/paulsen/})}

\begin{abstract}
We develop a symbol calculus for normal bimodule maps over a masa
that is the natural analogue of the Schur product theory. Using
this calculus we are able to easily give a complete description
of
the ranges of contractive normal bimodule idempotents that avoids
the theory of J*-algebras. 
 We prove that if $P$ is a normal
bimodule idempotent and $\|P\| < 2/\sqrt{3}$ then $P$ is a
contraction. We finish with some attempts at extending the symbol
calculus to non-normal maps.
\end{abstract}

\maketitle

\section{Introduction}

Let $\mathcal{H}$ be a Hilbert space and let $B(\cl H)$ denote
the
algebra of bounded linear operators on $\cl H$. If $P:B(\cl H)
\rightarrow B(\cl H)$ is a completely contractive idempotent,
then
it is known that its range $\cl M = P(B(\cl H))$ can be endowed
with a triple product that makes it a ternary ring of operators
(TRO). That is, there exists a unital $C^*$-algebra, $\cl A$, a
self-adjoint idempotent $q \in \cl A$, and a complete onto
isometry $\phi: \cl M \rightarrow q \cl A(1-q)$, such that the
triple product on $\cl M$ is given by $\{m_1,m_2,m_3 \} =
\phi(m_1)\phi(m_2)^*\phi(m_3).$

The presence of the map $\phi$ means that the triple product, in
general, has
little relation to the product inherited from $B(\cl H)$.
However, when $P$ is also assumed to be a bimodule map with
respect to a masa in
$B(\cl H)$, then Solel \cite{sol} has shown that considerably
more is true.

Let $\cl D\subseteq B(\cl H)$ be a masa acting on a separable
Hilbert space $\mathcal{H}$ and let $P:B(\cl H)\rightarrow B(\cl
H)$ be a contractive idempotent $\mathcal{D}$-bimodule map. In
case $P$ is also w*-continuous, Solel \cite{sol} proves that the
range $\mathcal{M}=P(B(\mathcal{H))}$ is actually a ternary ring
of operators (TRO) in the product inherited from $B(\cl H)$,
i.e.,
it satisfies $\mathcal{MM^*M\subseteq M}$. We summarize this
stronger condition by saying that $\cl M$ is a {\bf sub-TRO} of
$B(\cl H)$.

We give a new proof of Solel's result, that does not need any
J*-algebra 
theory, and in fact obtain the stronger conclusion that
$\mathcal{M}$ must be the w*-closed sum of ``full corners'', $\cl
M =\sum_n B(\cl H_n,\cl K_n)$ where $\{\cl H_n\}$ and $\{\cl
K_n\}$ are pairwise orthogonal subspaces of $\cl H$ and $\cl K$
respectively. From this description of $\cl M$ it is immediate
that $\cl M \cl M^* \cl M \subseteq \cl M$.

In addition, representing $\mathcal{D}$ as the multiplication
masa
of a standard Borel space $(X,\mu)$, we show that the
$\omega$-support (see \cite{eks}) of the range $\mathcal{M}$ of a
w*-continuous idempotent $\cl D$-bimodule map (whether
contractive
or not) is $\omega$-open (as well as $\omega$-closed). It follows
\cite{eks} that the reflexive cover $\mathrm{Ref} (\cl M)$ of
$\cl
M$ is strongly reflexive.

We begin by examining the case when $\cl H= \ell^2$ and $\cl D=
\ell^{\infty}$ in some detail. In this case every $\cl
D$-bimodule
map $\Phi$ is well-known to be given by Schur multiplication
against a fixed matrix $A=(a_{i,j})$, that is, $\Phi(X) = A*X=
(a_{i,j}x_{i,j}).$ We denote this map by $\Phi = S_A$.

Note that $\Phi \circ \Phi = \Phi$ is clearly equivalent to
$a_{i,j}^2 = a_{i,j}$ and hence each entry of $A$ must be either
0
or 1. Thus every idempotent $\Phi$ can be identified in a
one-to-one fashion with a subset $E \subseteq \bb N \times \bb N$
where $E= \{ (i,j) : a_{i,j}=1 \}$ and we write $A = \chi_E$,
where we regard the matrix $A$ as a function of two variables.

The problem of determining ranges of bounded bimodule
projections,
becomes one of determining which subsets of $\bb N \times \bb N$
give rise to bounded
bimodule projections.

When a non-zero idempotent Schur product map has norm less than
$2/\sqrt{3},$ then Livshits \cite{liv} proves that it must
actually be of norm 1.  We present a somewhat different proof of
this result, based on a combinatorial ``3 of 4'' lemma, that gives
a description of these sets that generalizes more readily from the
case of discrete masa's to continuous masa's.

Very little is known about the structure of sets such that $\Phi$
is only a bounded projection, but we give, hopefully, a little
insight into this problem.

In the third section, we use Haagerup's results \cite{unhaag} to
develop a symbol calculus for w*-continuous bimodule maps over
general masa's. One of the main advantages of our approach is that
the symbol calculus allows proofs given in the discrete case to
carry over to arbitrary masa's. In particular, we prove that every
non-zero w*-continuous bimodule projection of norm less than
$2/\sqrt{3}$ is actually of norm 1. Thus, the set of possible
norms of w*-continuous bimodule projections is not a connected
subset of the reals. It is not known if the same phenomenon holds
for all bounded bimodule projections.

In the final section we study masa bimodule maps that are not
w*-continuous. Solel \cite{sol} conjectures that the range of an
idempotent contractive masa bimodule map $P$ will still be a
sub-TRO even when $P$ is not w*-continuous. We use the actions of
such maps on normalizers of the bimodule as an approach to the
problem of extending our symbol calculus beyond the
w*-continuous case.

If $P$ is a contractive idempotent (hence $\Vert
P\Vert=1$) and its range is a $\cl D$-bimodule, then $P$ is
automatically a $\cl D$-bimodule map (\cite{sol}).

The authors would like to thank Ken Davidson, Allan Donsig, Gilles
Pisier and Ivan Todorov for various observations that have
improved our exposition and results.

\section{The Discrete Case}

In this section we develop the case where the masa is totally
atomic, so our Hilbert space may be represented as $\ell^2$ and
the masa as $\ell^{\infty}$ acting in the usual fashion as
diagonal matrices. Identifying $\ell^2 = L^2(\bb N, \mu)$, leads
to the identification of $\ell^{\infty} = L^{\infty}(\bb N, \mu)$
acting as multiplication operators on this space of functions.

\begin{defn}
Let $X$ and $Y$ be sets and let $E \subseteq X \times Y$. We say
that $E$ has the {\bf 3 of 4 property} provided that given any
distinct pair of points $x_1 \ne x_2$ in $X$ and any pair of
distinct points $y_1 \ne y_2$ in $Y$, whenever 3 of the 4 ordered
pairs $(x_i,y_j)$ belong to $E$ then the fourth ordered pair
belongs to $E$ also.
\end{defn}

\begin{lemma}\label{3of4}
Let $X$ and $Y$ be sets and let $E \subseteq X \times Y$. If $E$
has the 3 of 4 property, then there exists an index set $T$,
disjoint subsets $\{ X_t \}_{t \in T}$ of $X$, and disjoint
subsets $\{ Y_t \}_{t \in T}$ of $Y$ such that
$$E= \cup_{t \in T} X_t \times Y_t .$$
\end{lemma}
\begin{proof}
Define a relation on $X$ by $x_1 R x_2$ if and only if there
exists $y \in Y$ such that $(x_1,y)$ and $(x_2,y)$ are both in
$E.$ The 3 of 4 property ensures that $R$ is transitive and hence
is an equivalence relation on the subset $X_0= \{ x \in X: (x,y)
\in E \text{ for some } y\in Y \}$ of $X$. Let $\{ X_t \}_{t \in
T}$ denote the collection of equivalence classes of $X_0$.

Define $Y_t= \{ y \in Y: (x,y) \in E \text{ for some } x \in X_t
\}$. Again the 3 of 4 property ensures that the sets $Y_t$ are
disjoint with union equal to $$Y_0= \{ y \in Y: (x,y) \in E
\text{
for some } x\in X \}.$$

One final use of the 3 of 4 property shows that $E= \cup_{t \in
T}
X_t \times Y_t .$
\end{proof}

The following result is also Lemma 3 of \cite{liv},  but the proof
is short, so we include it.

\begin{lemma}\label{2byroot3}
Let $A = \begin{pmatrix} 1 & 1\\0 & 1 \end{pmatrix}$ and let $S_A
: M_2 \rightarrow M_2$ denote the map given as Schur product by
$A$. Then $\Vert S_A \Vert = 2/\sqrt{3}.$
\end{lemma}
\begin{proof}
Since the unitaries are the extreme points of the unit ball of
$M_2$, it is easily seen that $\Vert
S_{A}\Vert=\sup_{\theta}\left\|
\begin{pmatrix}
\cos\theta & \sin\theta\\
0 & \cos\theta
\end{pmatrix}
\right\|  $ and the result follows by computing this supremum.
\end{proof}

\begin{theorem}\label{discr}
Let $P: B(\ell^2) \rightarrow B(\ell^2)$ be a non-zero
$\ell^{\infty}$-bimodule map that is idempotent, let $\cl M =
P(B(\ell^2))$ denote the range of $P$ and assume that $\Vert P
\Vert < 2/\sqrt{3}.$ Then:
\begin{enumerate}
\item[(i)] $P=S_{\chi_E}$ where $E = \cup I_m \times J_m$, with
$\{ I_m \}$ and $\{ J_m \}$ countable collections of disjoint
subsets of $\bb N$,
\item[(ii)] $\cl M = \sum_m \chi_{I_m}B(\ell^2) \chi_{J_m}$ and
$\cl M \cl M^* \cl M \subseteq \cl M$,
\item[(iii)] $\Vert P \Vert = 1.$
\end{enumerate}
\end{theorem}

\begin{proof}
The fact that $P=S_{\chi_E}$ for some set $E \subseteq \bb N
\times \bb N$ was noted in the introduction. Choose any $x_1 \ne
x_2$ and $y_1 \ne y_2$ in $\bb N$ and consider the compression of
$P$ as a map from the span of $\{ e_{x_1}, e_{x_2} \}$ to the
span
of $\{ e_{y_1}, e_{y_2} \}$. Since $\Vert P \Vert < 2/\sqrt{3}$,
by the above lemma $E$ will have the 3 of 4 property and hence by
the first lemma be of the form given in (i). It is now obvious
that $\cl M$ will have the form claimed in (ii). The second
assertion in (ii) is immediate from this.

Alternatively, to see the second assertion in (ii), note that it
is enough to assume that the matrix units $E_{i,j}, E_{k,l}$ and
$E_{m,n}$ are in $\cl M$ and prove that $E_{i,j}E_{k,l}^*E_{m,n}$
is in $\cl M$. But this product will be 0 unless $j=l$ and $k=m$
in which case the product is $E_{i,n}$. However, in this case we
have that $(i,j), (k,j),(k,n)$ belong to $E$ and so again by the
3
of 4 property $(i,n)$ is in $E$ and so $E_{i,n} \in \cl M$.

Finally, to prove (iii), let $\{e_m \}$ denote the usual basis of
$\ell^2$, set $x_i = e_m$ when $i \in I_m$, set $y_j = e_m$ when
$j \in J_m$ and note that $\chi_E(i,j)= \langle x_i,y_j \rangle.$
Thus, $\Vert P \Vert \le 1$ by the theorem characterizing the
norms of Schur product maps, see for example \cite{pa}.
\end{proof}

Theorem 4(iii), with a somewhat longer proof, is contained
explicitly in \cite{liv} and (i) and (ii) can be obtained from the
proof given there.

\medskip

It is well known that the range of a completely contractive
projection is completely isometrically isomorphic to a TRO on
some
Hilbert space $\cl H$. This result induces a triple product on
the
range, but it is generally not the triple product given by the
original representation of the range as a subspace of $B (\cl
H)$.
Note that by (ii), we have that the range of $P$ is a TRO {\it in
the original triple product}, i.e., that the range is a sub-TRO
of
$B(\ell^2)$.

\begin{remark}
By the above result we see that the set of possible norms of
bimodule projections does not contain the interval from 1 to
$2/\sqrt{3}$. This makes the structure of this set somewhat
intriguing. Davidson has observed that the set of possible norms
is closed under product and under the taking of suprema. By a
result of Bhatia, Choi and Davis \cite{ch}, the number 2 is one
of the limit points
of this set. Other than these facts, not much seems to be known
about this set.
\end{remark}

If we let $\Delta = \{ E \subseteq \bb N \times \bb N: S_{\chi_E}
\text{ is bounded } \}$, then, by the results characterizing the
norms of Schur product maps, it is easily seen that $E \in
\Delta$
if and only if there exist bounded sequences of vectors, $\{ x_i
\}$ and $\{ y_j \}$, such that $\chi_E(i,j) = \langle x_i,y_j
\rangle$, but this characterization seems to be of little help in
obtaining other conditions that characterize the sets in $\Delta
.$

It is not hard to show that $\Delta$ is an algebra of sets and so
contains the algebra generated by the sets given by the above
theorem. It is not currently known whether $\Delta$ equals the
latter algebra.

\section{A Functional Calculus}

In the discrete case, every bounded bimodule map is given as a
Schur product map and so is automatically w*-continuous, but
this is not the case in general. In this section we develop a
functional calculus in the non-discrete case for w*-continuous
bimodule maps that allows us to treat these exactly like Schur
product maps and consequently obtain exact analogues of the
results of the previous section.

This functional calculus is different
from the one considered by Peller \cite{pe}, and appears to have
recently been discovered independently by Shulman and Kissin
\cite{KS}.

Let $\cl D \subseteq B(\cl H)$ be a masa acting on a separable
Hilbert space $\mathcal{H}$ and let $\Phi:B(\cl H)\rightarrow
B(\cl H)$ be a w*-continuous $\mathcal{D}$-bimodule map. Since
$\Phi$ is a bounded $\mathcal{D}$-bimodule map, a result of R.
Smith \cite{smi91} (see also \cite{davpo}) shows that $\Phi$ must
be completely bounded, and in fact
$\Vert\Phi\Vert_{cb}=\Vert\Phi\Vert$.

Haagerup \cite{unhaag} has shown that a w*-continuous completely
bounded $\mathcal{D}$-bimodule map such as $\Phi$ must be of the
form
\[
\Phi(T)=\sum_{n=1}^{\infty}F_{n}TG_{n}\quad(T\in B(\cl H))%
\]
for suitable $F_{n},G_{n}\in\cl D$ satisfying $\|\sum F_n F_n^*
\|<\infty$ and $\|\sum G_n^* G_n \|<\infty$.

Represent $\mathcal{D}$ as the multiplication masa of a standard
(finite) Borel space $(X,\mu)$ acting on
$\mathcal{H}=L^{2}(X,\mu)$. A standard null-set argument shows
that we may choose two families $\{ f_n \}, \{ g_n\}$ of
\emph{Borel functions } with $F_n=M_{f_n}$ and  $G_n=M_{g_n}$ for
each $n$, and such that the series $\sum|f_{n}(t)|^{2}$ and
$\sum|g_{n}(t)|^{2}$ converge for \emph{all } $t\in X$ boundedly
and in $L^{2}$ norm.

It follows that the series
\[
\phi(s,t)=\sum_{n}f_{n}(s)g_{n}(t)
\]
converges pointwise \emph{everywhere} to a Borel function.

\medskip

Conversely, let $f=(f_{1},f_{2},\ldots)$ and
$g=(g_{1},g_{2},\ldots)$ be (essentially) \emph{bounded} weakly
Borel measurable functions from $X$ into $\ell^{2}$. Since
$\ell^2$ is separable, weak Borel measurability and strong Borel
measurability are equivalent. Thus (i) each $f_{n}$ is an
essentially bounded complex-valued function and (ii)
\[
\sup_{s\in X}\Vert f(s)\Vert_{2}^{2}=\sup_{s\in
X}\sum_{n}|f_{n}(s)|^{2}\equiv B_{f}<\infty
\]
and $B_{g}\equiv\sup_{s\in X}\Vert g(s)\Vert_{2}^{2}<\infty$. It
follows that
\[
\Vert f\Vert^{2}\equiv\int\Vert
f(s)\Vert_{2}^{2}d\mu(s)=\sum_{n}\int
|f_{n}(s)|^{2}d\mu(s)<\infty
\]
and similarly for $g$; hence the function
\[
\phi(s,t)=\langle f(s),\bar{g}(t)\rangle=\sum_{n}f_{n}(s)g_{n}(t)
\]
defines an element of the projective tensor product $L^{2}(X,\mu)
\widehat{\otimes}L^{2}(X,\mu)$. Note that the series converges
pointwise absolutely and boundedly and also in the projective
norm. Thus the function $\phi$ is Borel on $X\times X$ and
(essentially) bounded.

Denoting by $F_{n}$ (resp. $G_{n}$) the multiplication operator
$M_{f_n}$ (resp. $M_{g_n}$) acting on $\cl H =L^2(X,\mu)$ we
observe that for every $T\in B(\cl H)$ the series
$\sum_{n}F_{n}TG_{n}$ converges in the w*-topology. This is well
known, and we sketch the proof for completeness.

The fact that the bounds $B_f$ and $B_g$ defined above are finite
shows that the infinite matrices $[G_1,G_2, \ldots ]^t$ and
$[F_1,F_2, \ldots ]$ define bounded operators $\cl H \to (\cl
H)^{(\infty)}$ and $(\cl H)^{(\infty)} \to \cl H$ with norms
$\sqrt{B_g}$ and $\sqrt{B_f}$ respectively. Thus the matrix
product $[F_1,F_2, \ldots ]T^{(\infty)}[G_1,G_2, \ldots ]^t$,
which we denote by $\Phi_{\phi}(T)$, defines an operator on $\cl
H$ whose norm is at most $\|T\| \sqrt{B_gB_f}$. This operator is
of course $\sum_{n}F_{n}TG_{n}$, the strong limit of its partial
sums $\Phi_{N}(T) =\sum_{n=1}^{N}F_{n}TG_{n}$; but since these
partial sums are uniformly bounded by $\|T\| \sqrt{B_gB_f}$, the
convergence is actually ultrastrong.

Clearly $\Phi_\phi$ is a $\mathcal{D}$-bimodule map. We claim
that
it is w*-continuous. For this, it suffices to show that it is
weak
operator continuous on the unit ball $B(\cl H)_1$ of $B(\cl H)$.
But if $T$ is a contraction, then for all
$\xi,\eta\in\mathcal{H}$, we have
\begin{align*}
& |\langle\Phi_{\phi}(T)\xi,\eta \rangle -\langle
\Phi_N(T)\xi,\eta\rangle|^2 = \left| \sum
_{n=N+1}^{\infty}\langle
TG_n\xi,F_{n}^{\ast}\eta\rangle\right|^2\\
&  \leq\left(  \sum_{n=N+1}^{\infty}\Vert
TG_{n}\xi\Vert^{2}\right) \left( \sum_{n=N+1}^{\infty}\Vert
F_{n}^{\ast}\eta\Vert^{2}\right) \leq\Vert T\Vert^2\left(
\sum_{n=N+1}^{\infty}\Vert G_n \xi\Vert^2 \right)B_f \|\eta\|^2
\\
&  \leq\left( \sum_{n=N+1}^{\infty}\int |g_{n}(s)|^{2}|\xi
(s)|^{2}d\mu(s)\right)  B_{f}\Vert\eta\Vert^2 \, .
\end{align*}
By dominated convergence, this tends to $0$ as $N\to \infty$. It
follows that the function
$T\rightarrow\langle\Phi_{\phi}(T)\xi,\eta \rangle$ is the
\emph{uniform} limit on $B(\cl H)_{1}$ of the weak operator
continuous functions
$T\rightarrow\langle\Phi_{N}(T)\xi,\eta\rangle$ and so is itself
weak operator continuous.

Note that the map $\Phi_\phi$ acts as a multiplication operator
on
kernels of Hilbert Schmidt operators:

\begin{prop}\label{hs}
Let $\phi(s,t)=\langle f(s),\bar{g}(t)\rangle$ where $f$ and $g$
are (essentially) bounded (weakly) Borel measurable functions
from
$X$ into $\ell^{2}$, and let $\Phi_{\phi}:B(\cl H)\rightarrow
B(\cl H)$ be as above. The map $\Phi_\phi$ leaves the space of
Hilbert Schmidt operators invariant. If $T\in B(\cl H)$ is a
Hilbert Schmidt operator with kernel $k\in L^2(X\times X)$, then
$\Phi_\phi(T)$ has kernel $M_{\phi}(k)=\phi k$. Thus $\Phi_\phi$
acts on $L^{2}(X\times X)$ as multiplication by $\phi$.
\end{prop}

\begin{proof} Let $T=T_k$ be Hilbert Schmidt operator with
kernel $k$. For $\xi,\eta\in\mathcal{H}$ we have \vspace{-1ex}
\begin{align*}
\langle \Phi_\phi(T_{k})\xi,\eta\rangle & =
\langle\sum_{n=1}^{\infty}F_nT_k G_n \xi,\eta\rangle
=\sum_{n=1}^{\infty}\langle
T_{k}G_n\xi,F_n^\ast \eta\rangle\\
&  =\sum_{n=1}^{\infty}\iint
k(x,y)g_n(y)\xi(y)f_n(x)\overline{\eta(x)}
d\mu(y)d\mu(x)\\
&
=\iint\sum_{n=1}^{\infty}f_n(x)g_n(y)k(x,y)\xi(y)\overline{\eta(x)}
d\mu(y)d\mu(x)\\
&  =\iint\phi(x,y)k(x,y)\xi(y)\overline{\eta(x)}d\mu(y)d\mu(x)\\
&  =\langle T_{\phi k}\xi,\eta\rangle
\end{align*}
and so $\Phi_\phi(T_{k})=T_{\phi k}$.
\end{proof}

The next result shows that a w*-continuous bimodule map $\Phi =
\Phi_\phi$ determines the function $\phi$ essentially uniquely.
Recall that a set $R \subseteq X \times X$ is called {\bf
marginally null} provided that there is a null set $N$ such that
$R \subseteq (N\times X)\cup (X\times N)$. We say that two
functions $\phi$ and $\psi$ on $X\times X$ are equal {\bf
marginally almost everywhere} and write $\phi = \psi$  m.a.e.
provided that the set of points where they are not equal is
marginally null.

\begin{theorem}\label{mae}
Let $\phi(s,t)=\langle f(s),\bar{g}(t)\rangle$ where $f$ and $g$
are (essentially) bounded (weakly) Borel measurable functions
from
$X$ into $\ell^{2}$, and let $\Phi_{\phi}:B(\cl H)\rightarrow
B(\cl H)$ be as above. The following are equivalent:

\begin{enumerate}
\item $\phi=0$ m.a.e.

\item $\phi=0$ a.e.

\item $\Phi_{\phi}=0$.
\end{enumerate}
\end{theorem}

\begin{proof} If the set \vspace{-1ex}
\[
R=\{(s,t)\in X\times X:\phi(s,t)\neq0\}
\]
is contained in a set of the form $N\times X\cup X\times N$,
where
$N\subseteq X$ is null, then of course the product measure of $R$
is $0$. Thus (1) implies (2).

To show that (2) implies (3), observe that if $T=T_k$ is a
Hilbert-Schmidt operator with (square-integrable) kernel $k$,
then
by Proposition \ref{hs} $\Phi_{\phi}(T_{k})=T_{\phi k}$.

It follows that if $\phi=0$ a.e. then $\Phi_{\phi}(T_k)=0$ for
any
Hilbert-Schmidt operator $T_k$. Since $\Phi_{\phi}$ is
w*-continuous, we obtain $\Phi_{\phi}=0$. Conversely if
$\Phi_{\phi}=0$ then $\phi=0$ a.e.

It remains to prove that if the set $R$ is null, then it must be
marginally null. For this, first observe that $R$ is (marginally
equivalent to) a countable union of Borel rectangles. We use an
argument of Arveson \cite{arv74}: The set
\[
\{(\xi,\eta)\in\ell^{2}\times\ell^{2}:\left\langle
\xi,\eta\right\rangle
\neq0\}
\]
is open in $\ell^{2}\times\ell^{2}$, and hence is a countable
union
$\cup
_{n}U_{n}\times V_{n}$ of open rectangles. Letting $A_{n}=\{s\in
X:f(s)\in
U_{n}\}$ and $B_{n}=\{t\in X:g(t)\in V_{n}\}$ we see that, since
$f,g:X\rightarrow\ell^{2}$ are Borel functions, the sets $A_{n}$
and
$B_{n}$
are Borel and
\[
R=\{(s,t):\left\langle f(s),g(t)\right\rangle
\neq0\}=\bigcup_{n}%
\{(s,t):\left\langle f(s),g(t)\right\rangle \in U_{n}\times
V_{n}%
\}=\bigcup_{n}A_{n}\times B_{n}%
\]
as claimed. Thus if the product measure of $R$ is $0$ we must
have
$\mu (A_{n})\mu(B_{n})=0$ for all $n\in\mathbb{N}$. If
$N_1=\cup\{A_n:\mu (B_n)\neq0\}$ and
$N_2=\cup\{B_n:\mu(A_n)\neq0\}$ then $\mu(N_1)=\mu(N_2)=0$ and
\[
R\subseteq N_{1}\times X\cup X\times N_{2}%
\]
which completes the proof.
\end{proof}
\begin{defn}
We let {\bf $NCB_{\cl D}(B(\cl H))$ } denote the algebra of
w*-continuous $\cl D$-bimodule maps from $B(\cl H)$ into
itself. Given a w*-continuous $\cl D$-bimodule map $\Phi$ as
above we call the m.a.e. equivalence class of the function $\phi$
obtained above the {\bf symbol} of $\Phi$ and denote it by
$\Gamma(\Phi).$
\end{defn}

\begin{corollary}
Let $\cl D$ be represented as the multiplication masa on a
standard Borel space $(X, \mu)$, let $\cl B_{mae}(X \times X)$
denote the algebra of bounded Borel functions on $X \times X$
modulo the marginally null functions. Then the map $\Gamma :
NCB_{\cl D}(B(\cl H)) \rightarrow \cl B_{mae}(X \times X)$ is a
one-to-one homomorphism onto the subalgebra of functions that can
be represented in the form $\phi(s,t) = \langle f(s),g(t)
\rangle$
for any bounded Borel measurable functions $f,g$ from $X$ into a
separable Hilbert space.
\end{corollary}

We call the map $\Gamma$ the {\bf functional calculus} for
w*-continuous $\cl D$-bimodule maps.

Armed with the functional calculus, we can readily generalize the
theorem of the previous section.

\begin{theorem}\label{cts}
Let $P: B(\cl H) \rightarrow B(\cl H)$ be a $\cl D$-bimodule map
that is idempotent and w*-continuous, let $\cl M = P(B(\cl
H))$
denote the range of $P$ and assume that $\Vert P \Vert <
2/\sqrt{3}.$ Then:
\begin{enumerate}
\item $\Gamma(P) =\chi_E$ where $E = \cup I_m \times J_m$, with
$\{ I_m \}$ and $\{ J_m \}$  countable collections of disjoint
Borel subsets of $X$,

\item $\cl M = \sum_m \chi_{I_m}B(L^2) \chi_{J_m}$ and $\cl M
\cl M^* \cl M \subseteq \cl M$,
\item $\Vert P \Vert = 1.$
\end{enumerate}
\end{theorem}
\begin{proof}
Let $\Gamma(P)= \phi$, since $P \circ P = P$, by the functional
calculus, $\phi^2 = \phi$ marginally almost everywhere. Thus, we
can pick a Borel subset $X_1$ of $X$ with $\mu(X \cap X_1^c) =0$
such that $\phi^2 = \phi$ on $X_1 \times X_1$. Hence, there is a
Borel subset $E$ of $X_1 \times X_1$ such that $\phi = \chi_E$ as
functions on $X_1 \times X_1$.

Thus, we may write $\chi_E(s,t) = \langle f(s),g(t) \rangle$
where
$f,g$ are functions into a separable Hilbert space with $\Vert
f(s) \Vert \Vert g(t) \Vert < 2/\sqrt{3}$ for all $s,t.$

By Lemma \ref{2byroot3}, if for any $s_1 \ne s_2$ and $t_1 \ne
t_2$, we have that 3 of the 4 values $\langle f(s_i),g(t_j)
\rangle$ are 1, then the fourth value must also be 1.

Hence the set $E$ satisfies the 3 of 4 property and so it must be
a union of disjoint rectangles as in Lemma \ref{3of4}. Say, $E =
\cup_{t \in T} I_t \times J_t.$

It remains to be shown that the indexing set $T$ for the union is
only countable and that each of the sets $I_t$ and $J_t$ are
Borel. As in the proof of Theorem \ref{mae}, the set
$E=\{(s,t)\in
X_1\times X_1 : \phi(s,t) \ne 0 \}$ can be written as a countable
union of Borel rectangles, say $E= \cup_n A_n \times B_n $.
Again,
by the equivalence relation used to define the sets $I_t$ and
$J_t$, if any point in a rectangle $A_n \times B_n$ is contained
in $I_t \times J_t$ then $A_n \times B_n \subseteq I_t \times
J_t.$ Hence, each set $I_t \times J_t$ is the union of the at
most
countably many Borel rectangles that are contained in it and
consequently is itself a Borel rectangle. Moreover, the set $T$
can be placed in a one-to-one correspondence with a partition of
the integers, and hence is countable.

The remainder of the proof proceeds as in the proof of Theorem
\ref{discr}.
\end{proof}

As in the discrete case we have that $\cl M$ is a sub-TRO of
$B(L^2)$, a result obtained by Solel \cite{sol}.

Just as in the discrete case, very little is known about bimodule
projections of greater norm.

More importantly, very little is known about contractive bimodule
projections that are not w*-continuous. Such projections do
exist, for example projections onto the masa $\cl D$ exist and
when $\cl D$ is not discrete, these cannot be 
not
w*-continuous. Solel \cite{sol} conjectures that the range
$\cl
M$ of {\it any} contractive $\cl D$-bimodule projection satisfies
$\cl M \cl M^* \cl M \subseteq \cl M.$

\bigskip

We now turn our attention to further properties of the symbol
calculus and of bounded w*-continuous idempotents.

Recall \cite{eks} that a subset $E \subseteq X\times X$ is said
to
be $\omega$-\emph{open } if it differs from a countable union of
Borel rectangles by a marginally null set, and is
$\omega$-\emph{closed } if its complement is $\omega$-open. Thus
the set $E$ in Theorem \ref{cts} is $\omega$-open. The following
Proposition strengthens this result and provides an alternative
approach:

\begin{prop}
Let $P \in NCB_{\cl D}(B(\cl H))$ be an idempotent with symbol
$\Gamma(P) = \chi.$ Then there exists an $\omega$-open and
$\omega$-closed set $A\subseteq X\times X$ such that
$\chi=\chi_{A}$ marginally almost everywhere.
\end{prop}

\begin{proof} Notice first that, in the terminology of
\cite{eks}, any element $\phi$ of the projective tensor product
$L^2(X)\widehat{\otimes} L^2(X)$ is $\omega$-\emph{continuous},
that is, $\phi^{-1}(U)$ is $\omega$-open in $X\times X$ for any
open set $U \subseteq \bb C$ \cite[Theorem 6.5]{eks}.

Since $P$ is idempotent, so is its induced operator $M_{\chi}$ on
$L^{2}(X\times X)$ (see Proposition \ref{hs}). It follows that
$\chi^{2}=\chi$ almost everywhere, i.e. the set
\[
B=\{(s,t):\chi^{2}(s,t)-\chi(s,t)\neq0\}
\]
has product measure zero. On the other hand, since $\chi \in
L^2(X)\widehat{\otimes} L^2(X)$, the function $\chi$ is
$\omega$-continuous hence so is $\chi^2 -\chi$. Thus $B$ must be
$\omega$-open, in other words marginally equivalent to a
countable
union of rectangles. The fact that $B$ has product measure zero
now implies, as noted earlier, that it is actually marginally
null.  Replacing $X$ by a suitable Borel subset $X_1$ such that
$\mu (X\cap X_1^c)=0$, we may assume that $B=\emptyset$, i.e.
that
$\chi^{2}(s,t)=\chi(s,t)$ for \emph{all} $(s,t)\in X\times X$.
Thus letting $A=\chi^{-1}(\{1\})$ we see that $A$ is
$\omega$-closed (since $\chi$ is $\omega$-continuous); but
$A^{c}=\chi^{-1}(\{0\})$ is also $\omega$-closed.
\end{proof}

It is shown in \cite{eks} that, given any space $\cl S$ of
operators on $L^2(X)$, there exists an $\omega$-closed set
$\Omega$, minimal up to marginally null sets, that
\emph{supports}
all elements of $\cl S$, in the sense that if a Borel rectangle
$\alpha \times\beta$ doesn't meet $\Omega$ then $M_{\beta }\cl S
M_{\alpha}=\{0\}$ (here $M_{\beta}\in\mathcal{B}(L^{2}(X))$
denotes the projection onto $L^{2}(\beta))$. This set is called
the $\omega$-\emph{support } of $\cl S$.

\begin{prop}
Let $P \in NCB_{\cl D}(B(\cl H))$ be an idempotent with symbol
$\Gamma(P) = \chi_A.$ Then the set $A$ is (marginally equivalent
to) the $\omega$-support of $\cl M = P(B(\cl H))$.
\end{prop}

\begin{proof} It is to be shown that a Borel rectangle $\alpha
\times\beta$ has marginally null intersection with $A$ if and
only
if $M_\beta \cl M M_\alpha =\{0\}$. Note that the relation
$M_{\beta }\mathcal{M}M_{\alpha}=\{0\}$ is equivalent to
$M_{\beta}P(T)M_{\alpha}=0$ for all $T\in B(\cl H)$. But, since
the map $T\rightarrow M_{\beta}P(T)M_{\alpha}$ is w*-continuous,
this is equivalent to $M_{\beta}P(T)M_{\alpha}=0$ for all
\emph{Hilbert Schmidt }$T=T_{k}$. By Proposition \ref{hs}
$M_{\beta}P(T_{k})M_{\alpha}=M_{\beta }T_{\chi
k}M_{\alpha}=T_{h}$, where $h=\chi_{\alpha\times\beta}\chi_{A}k$.
Thus the relation $M_{\beta}P(T)M_{\alpha}=0$ holds for all
Hilbert Schmidt $T=T_{k}$ if and only if the set
$(\alpha\times\beta)\cap A$ has product measure zero. But since
this set is $\omega$-open, as shown in the proof of the last
proposition this can only happen when $(\alpha\times\beta)\cap A$
is \emph{marginally null}.
\end{proof}

Since $A$ is $\omega $-open, it follows \cite[Theorem 6.11]{eks}
that the reflexive cover $\mathrm{Ref}(\mathcal{M})$ is in fact
strongly reflexive, and is the strong closure of the linear span
of the finite rank operators supported in $A$. In case $P$ is
actually contractive, this of course follows immediately from the
fact that the range of $P$ is a direct sum of full corners
(Theorem \ref{cts}).

\section{Module maps and Normalizers}

We conjecture that it is possible to extend the functional
calculus of the previous section beyond the w*-continuous
module maps to a homomorphism from all of $CB_{\mathcal{D}}(B(\cl
H))$ into an algebra of functions. When $\mathcal{D}$ is a
continuous masa, then \cite{HW} proves that not only is
$CB_{\mathcal{D}}(B(\cl H))$ not abelian, but that its center is
exactly the w*-continuous maps, $NCB_{\cl D}(B(\cl H))$.

Since the range of the homomorphism that we hope to construct
would be contained in an abelian algebra, then necessarily this
extended functional calculus would need to annihilate the
commutator ideal. For this reason, we seek a representation of
this algebra that annihilates the commutator ideal and a general
understanding of the commutator ideal. We believe that the key to
finding such a representation is to understand the action of
module maps on normalizers.

Recall that an operator $T\in B(\cl H)$ is said to
\emph{normalize} a masa $\mathcal{D}\subseteq B(\cl H)$ if
$T^{\ast}\mathcal{D}T\subseteq \mathcal{D}$ and
$T\mathcal{D}T^{\ast}\subseteq\mathcal{D}$. If $T=V|T|$ is the
polar decomposition, then $T$ normalizes $\mathcal{D}$ if and
only
if $|T|\in\mathcal{D}$ and $V$ normalizes $\mathcal{D}$.

The set of all normalizers of $\mathcal{D}$ is a selfadjoint
semigroup (under composition) and it contains $\cl D$. Therefore
its norm-closed linear span is a unital C*-algebra, $\cl B$ say.
If we decompose $\cl D$ into its totally atomic and continuous
parts and let $\cl H = \cl H_a \oplus \cl H_c$ be the
corresponding orthogonal decomposition of the underlying Hilbert
space, then any partial isometry that normalizes $\cl D$ must
leave each of these subspaces invariant. It follows that $\cl B$
has at most one nontrivial invariant projection, corresponding to
this decomposition of $\cl D$. For a $\cl B$-projection must on
the one hand commute with $\cl D$ hence must belong to it; on the
other hand if a projection in $\cl D$ is invariant under all
normalizers, then it must be trivial or the projection onto one
of
these parts. Thus, if we assume that $\cl D$ is either atomic or
continuous, then the w*-closure of $\cl B$ is $B(\cl H)$. We will
show below that $\cl B$ is always proper.

When $\cl D$ is atomic, then it is clear that $\cl B$ contains
all
the corresponding matrix units and hence all of the compact
operators.

\begin{remark}
If $\cl D$ is continuous, does $\cl B$ contain compact operators?
Since $\cl B$ is irreducible, it follows that if it contains some
nonzero compact operator then it must contain them all; for the
set of all compact operators in $\cl B$ is then a nonzero ideal
of
the irreducible algebra $\cl B$, hence must itself be
irreducible.
\end{remark}

\begin{prop}
\label{one}Let $\Phi\in CB_{\cl D}=CB_{\cl D}(B(\cl H))$. For
each
normalizer $T$ of $\cl D$ there exists $D_{T}\in\cl D$ such that
\[
\Phi(T)=D_{T}T.
\]
In fact $D_{T}$ depends only on the partially isometric factor
$V$
in the polar decomposition $T=V|T|$ of $T$.
\end{prop}

Thus $\Phi(Vf(|T|))=D_{T}Vf(|T|)$ for any Borel function defined
on $\sigma(|T|)$.

\begin{proof}
Since $V$ itself normalizes the masa while $|T|$ belongs to it,
$\Phi(V|T|)=\Phi(V)|T|$ by
the module property of $\Phi$. Hence it suffices to deal with a
normalizing partial isometry $V$. Now for each $D\in\mathcal{D}$
the element
\[
\Delta(V)=V^{\ast}DV
\]
belongs to $\cl D$, and so
\begin{align*}
V^{\ast}D  &
=V^{\ast}(VV^{\ast})D=V^{\ast}DVV^{\ast}=\Delta(V)V^{\ast}\\
\text{and}\qquad V\Delta(V)  &  =(VV^{\ast})DV=DVV^{\ast}V=DV
\end{align*}
Therefore
\begin{align*}
\Phi(V)V^{\ast}D  &
=\Phi(V)\Delta(V)V^{\ast}=\Phi(V\Delta(V))V^{\ast}\\
&  =\Phi(DV)V^{\ast}=D\Phi(V)V^{\ast}.
\end{align*}
We have shown that $\Phi(V)V^{\ast}$ commutes with $\cl D$ and
hence there exists $D_{V}\in\cl D$ so that $ \Phi(V)V^*=D_V$
hence
\[
D_{V}V=\Phi(V)V^{\ast}V=\Phi(VV^{\ast}V)=\Phi(V)
\]
as required.
\end{proof}

The above result does provide some sort of ``symbol'' for $\Phi$.
If we assume that we are in the situation of the last section,
namely, that $\cl H = L^2(X, \mu)$, and that $\Phi$ is a
w*-continuous $L^\infty$-module map with symbol $\phi(s,t)$,
then given any function $h$ on $X$ that gives rise to a bounded
composition operator $C_h$ on $\cl H$, it is easy to see that
$\Phi(C_h) = M_k C_h$, where $k(s) = \phi(s,h(s))$ a.e. This
follows from the easily verified fact that $M_fC_hM_g = M_kC_h$
where $k(s)=f(s)g(h(s))$ a.e. and that $\Phi$ can be represented
as a sum of such elementary operators.

Thus, the function $k$ tells us the value of the symbol
restricted
to the graph of $h.$

By the above result, even in the case that $\Phi$ is not
w*-continuous, we see that for each bounded composition operator
it is still true that $\Phi(C_h) = M_k C_h$ for some function $k$
and the problem becomes, whether or not we can construct a symbol
$\phi(s,t)$ such that $k(s) = \phi(s, h(s))$ holds for every $h$
that gives rise to a bounded composition operator.

\medskip

For the following results, the reader should keep in mind that it
is quite likely that $\cl B = \cl B+ \cl K.$

\begin{prop}
Let $\cl B \subseteq B(\cl H)$ be the C*-algebra generated
by the normalizers of $\cl D$ and let $\cl K$ denote the
ideal of compact operators. Then completely bounded
$\cl D$-bimodule maps leave
$\cl B,\, \cl K$ and $\cl B+ \cl K$ invariant and the algebra
obtained by restricting the $\cl D$-bimodule maps to any
of these subspaces is commutative.
\end{prop}

\begin{proof} The fact that $\cl B$ is left invariant by every
$\cl D$-bimodule map follows from Proposition \ref{one}. Let
$\Phi_1,\Phi_2\in CB_{\cl D}(B(\cl H))$. If $T\in B(\cl H)$ is a
normalizer of $\cl D$, by Proposition \ref{one} there exist
$D_{1},D_{2}\in\cl D$ s.t.
\[
\Phi_{1}(T)=D_{1}T,\qquad\Phi_{2}(T)=D_{2}T.
\]
Since $\cl D$ is abelian, it follows that
\begin{align*}
\Phi_{1}(\Phi_{2}(T))  &
=\Phi_{1}(D_{2}T)=D_{2}\Phi_{1}(T)=D_{2}%
D_{1}T\\
&
=D_{1}\Phi_{2}(T)=\Phi_{2}(D_{1}T)=\Phi_{2}(\Phi_{1}(T)).\nonumber
\end{align*}
Thus the commutator $\Phi_{1}\circ\Phi_{2}-\Phi_{2}\circ\Phi_{1}$
vanishes on normalizers; hence it vanishes on $\cl B$ by
linearity
and continuity.

Hence, the algebra obtained by restricting $\cl D$-bimodule maps
to $\cl B$ is commutative.

Proposition \ref{hs} shows that every normal bimodule map
leaves $\cl K$ invariant. If we decompose an arbitrary bimodule
map into its normal and singular part, then the singular part
annihilates $\cl K$ and so also leaves it invariant. Thus, $\cl
K$ is left invariant by every $\cl D$-bimodule map and the
algebra obtained by restricting $\cl D$-bimodule maps to $\cl
K$ is exactly $NCB_{\cl D}(B(\cl H))$.

The result for $\cl B+ \cl K$ follows from the first two parts.
\end{proof}

\begin{corollary}
The restriction maps
\[
R_1:CB_{\cl D}(B(\cl H))\rightarrow
CB_{\cl D}(\cl B%
):\Phi\rightarrow\Phi|_{\cl B}%
\]
and
\[
R_2:CB_{\cl D}(B(\cl H))\rightarrow
CB_{\cl D}(\cl B+ \cl K
):\Phi\rightarrow\Phi|_{\cl B+ \cl K}%
\]
are onto, quotient maps that vanish on the commutator ideal of
$CB_{\cl D}(B(\cl H))$.
\end{corollary}
\begin{proof}
The fact that $R_1$ vanishes on commutators is clear. By
Wittstock's extension theorem for completely bounded bimodule
maps
\cite{W}, every map in $CB_{\cl D}(\cl B)$ has an extension to a
map of the same completely bounded norm in $CB_{\cl D}(B(\cl
H)).$
The fact that an extension exists means that $R_1$ is onto, and
the fact that the (cb) norm remains the same means that $R_1$ is
a
(complete) quotient map.

The proof for $\cl B+ \cl K$ is identical.
\end{proof}

\noindent\textbf{Conjecture.}\emph{ We conjecture that the
kernels
of $R_1$ and of $R_2$ are both equal to the (closed) commutator
ideal. }

The above result also provides some insight into the structure of
$\cl B.$

\begin{corollary}
When $\cl D$ is not totally atomic, the C*-algebras $\cl B$ and
$\cl B+ \cl K$ are proper (w*-dense) subalgebras of $B(\cl H)$.
\end{corollary}

\begin{proof} If $\cl B+ \cl K$ were
equal to $B(\cl H)$ then all $\cl D$-bimodule maps would commute.
But this is not the case when $\cl D$ is not totally atomic by
\cite{HW}. Hence $\cl B$ is also proper in this case.
\end{proof}

The fact that the set of $\cl D$-bimodule maps is non-commutative
when $\cl D$ is not totally atomic can also be seen directly.
For
example, if $\Phi_{1},\Phi_{2}$ are two (norm-continuous)
distinct
projections onto the masa itself, then their commutator is their
difference. To see that two such distinct projections exist it
suffices, by splitting $\cl D$ into its atomic and continuous
part, to assume $\cl D\simeq L^{\infty}[0,1]$ (see for example
\cite{tak}, 
Theorem III.1.22). For this case two distinct projections onto
$\cl D$ can be found for example in \cite{dav}, 
Algebras,
8.9 and 8.10.

\bigskip

The above proof does not apply to the totally atomic case, since
in this case the $\cl D$-bimodule maps do commute. Moreover, as
remarked above, $\cl B$ contains all compacts. However, another
argument shows that $\cl B$ is still a proper subalgebra of
$B(\cl H).$ We argue for the case that $\cl D = \ell^{\infty}.$

In this case, we have that $\cl B \subseteq B(\ell^2)$ is the
C*-algebra generated by the
diagonal $\cl D = \ell^\infty$ and its normalizers. It is
convenient to label the basis of $\ell^{2}$ as $e_{0},e_{1},
\ldots$. For each $A=[a_{i,j}] \in B(\ell^2)$, we denote by
$a_{n}
\, (n\in\bb Z)$ its $n$-th diagonal, i.e. $a_{n} =
(a_{i,j})_{i-j=n} \in \ell^{\infty}$, so that formally
$A=\sum_{n\in\bb Z} M_{a_{n}}S_{n}$ where $M_{a}\in B(\ell^{2})$
is the operator of multiplication by $a\in \ell^{\infty}$ and we
set $S_{n} =S^{n}$ for $n \ge0$ and $S_{n} = (S^{*})^{-n}$ for $n
\le0$, where $S$ is the unilateral shift.

\begin{lemma}
Let $m:\ell^{\infty}\rightarrow\bb C$ be a Banach limit, i.e. a
translation invariant state. If $A=\sum_{|n|\leq
N}M_{a_{n}}S_{n}$
has
finitely many nonzero diagonals, we set
$\Psi_{o}(A)=\sum_{|n|\leq N}%
m(a_{n})S_{n}$. Then $\Psi_{o}$ extends to a completely positive
projection $\Psi:B(\ell^2)\rightarrow B(\ell^2)$ onto the
Toeplitz
operators.
\end{lemma}

\begin{proof} For $A \in B(\cl H)$, consider the sesquilinear
form
\[
\cl H \times\cl H \to\bb C : (x,y) \to m( (\langle AS^{n} x ,
S^{n} y \rangle)_{n\in\bb N}).
\]
Since
\[
| m( (\langle AS^{n} x , S^{n} y \rangle)_{n\in\bb N}) | \le\|
m\|
\sup_{n\in\bb N} |\langle AS^{n} x , S^{n} y \rangle| \le\| A\|
\|x\|\|y\| ,
\]
(because $\| m \| =1$ and $\|S^{n}\|=1$) it is clear that this
form is bounded, so there exists a unique $\Psi(A) \in B(\cl H)$
such that
\[
\langle\Psi(A) x,y \rangle= m( (\langle AS^{n} x , S^{n} y
\rangle
)_{n\in\bb N}).
\]
Note that the $k$-th diagonal of $A$ is the sequence $a_{k} =
(\langle Ae_{k+j}, e_{j}\rangle)_{j\in\bb N}$. Now for all $i,k
\in\bb N$, the sequence $(\langle Ae_{k+i+n},
e_{i+n}\rangle)_{n\in\bb N}$ is a translate of the sequence
$a_{k}$, and so they have the same mean. Thus
\[
\langle\Psi(A)e_{k+i}, e_{i}\rangle= m( (\langle AS^{n} e_{k+i} ,
S^{n} e_{i} \rangle)_{n}) = m( (\langle Ae_{n+k+i} , e_{n+i}
\rangle)_{n}) = m(a_{k})
\]
hence the $k$-th diagonal of $\Psi(A)$ is the constant sequence
$(m(a_{k}), m(a_{k}), \ldots)$. Therefore $\Psi$ extends
$\Psi_{o}$.

It is clear that $\Psi$ is a linear projection onto the set of
Toeplitz operators. Finally, $\Psi$ is $k$-positive for all
$k\in\bb N$. Indeed if $[A_{ij}]\in M_{k}(B(\cl H))$ is positive
and $x=(x_{i})\in \cl H^{k}$, we have
\begin{align*}
&  \left\langle [\Psi(A_{ij})]x,x\right\rangle
=\sum_{ij}\left\langle
\Psi(A_{ij})x_{j},x_{i}\right\rangle =\sum_{ij}m((\langle
A_{ij}S^{n}%
x_{j},S^{n}x_{i}\rangle)_{n})\\
&  =m\left(  \left(  \sum_{ij}\langle
A_{ij}S^{n}x_{j},S^{n}x_{i}%
\rangle\right)  _{\!n}\right)  =m\left(  \left(  \langle\lbrack
A_{ij}%
]\tilde{S}^{n}x,\tilde{S}^{n}x\rangle\right)  _{n}\right)
\geq0\\
&
\end{align*}
(where
$\tilde{S}^{n}x=(S^{n}x_{1},S^{n}x_{2},\ldots,S^{n}x_{k})$).
\end{proof}

\begin{prop}
Let $\cl B$ and $\Psi$ be as above.
The image $\Psi (\cl B)$ contains only Toeplitz operators with
continuous symbol; hence $\cl B$ cannot be all of
$B(\ell^2)$.
\end{prop}

\begin{proof}
Given a  partially isometric normalizer $V$ of $\cl D$ let $
(v_n)$ be its sequence of diagonals. Note that since $V$ has at
most one 1 per row and column the (coordinatewise) sum
$\sum_{n\in
\bb Z} v_n$ is dominated by the vector $(2,2, \ldots) \in
\ell^\infty$. It follows that for each $N\in \bb N$,
$$0 \le \sum_{|n|\le N} m(v_n) = m \left(\sum_{|n|\le N} v_n
\right) \le m((2,2, \ldots))=2$$ and so $0\le \sum_{n\in \bb Z}
m(v_n) \le 2$. It follows that $\Psi (V) = T_f$ where $\sum
|\hat{f}(n)| \le 2$ and so $f$ is continuous (its Fourier series
converges absolutely).

Consequently for any diagonal operator $D$ of norm $M$ say, the
sum of the absolute values of the diagonals of $DV$ will be
dominated by the vector  $(2M,2M, \ldots)$ and so $\Psi (DV)=T_f$
where $\sum_n| \hat{f}(n)| \le 2M$. Hence $f$ will be a
continuous
function on the circle. By linearity and continuity of $\Psi$,
the
result follows.
\end{proof}

\end{document}